\documentclass[leqno,12pt]{article} 
\usepackage{amsmath, amssymb}
\usepackage{amsthm} 
\theoremstyle{plain} 
\newtheorem{theorem}{\indent\sc Theorem}[section]
\newtheorem{lemma}[theorem]{\indent\sc Lemma}

\theoremstyle{definition} 

\newcommand{\proj}{\mathbb P}
\newcommand{\reel}{\mathbb R}
\newcommand{\complex}{\mathbb C}

\newcommand{\é}{\'{e}}

\title{\uppercase{Special Lagrangian manifolds obtained from complex Grassmannians.}}
\author{\textsc{A. Ben Abdesselem and P. Cabau}}
\date{Inst. Math. Jussieu (UMR 7586)\\
Case 247 - 4, place Jussieu,
75252 Paris Cedex FRANCE.\\
E-mail address~: benabdes@math.jussieu.fr\\
Lab. Ing. Math.\\
Ecole Polytechnique de Tunisie, La Marsa 2070,
TUNISIA\\
E-mail address~: patrickcabau@yahoo.fr}

\begin{document}

\maketitle


\begin{abstract}
This paper gives an example of special Lagrangian manifold obtained
from a hypersurface of a complex Grassmannian with vanishing first
Chern class. The obtained manifold is a $1$-torus bundle over the
two dimensional real projective space. Such manifolds are
interesting for mirror symmetry theory. Other examples of the same
type are provided at the end of this article.

\end{abstract}

{\section*{Introduction.}}

In 1982, when F. Harvey and H. Lawson introduced Special Lagrangian
manifolds in \cite{HL}, their main interest was calibration
problems. Now, these manifolds give another approach to mirror
symmetry and string theory and become crucial in these fields. A
conjecture due to Strominger, Yau and Zaslow in \cite{SYZ} explains
mirror symmetry in a fairly mathematical way. An important class of
examples of special Lagrangian manifolds $L$ may be found among the
submanifolds of a $n$ complex dimensional manifold $M$, which is
compact k\ählerian with vanishing first Chern class. Then, the
definition is given by the following points:
\begin{itemize}
\item the symplectic form $\omega$ associated to the k\ählerian
structure restricted to the $n$ real dimensional submanifold $L$
must identically vanish, i.e. $i^{*}\omega=0$, where $i:L\rightarrow
M$ is the canonical injection. So the maximal isotropic submanifold
$L$, is endowed by a Lagrangian structure.
\item There exists a $(n,0)$-holomorphic volume form $\Omega$ such
that
$$ \Omega \wedge \overline{\Omega}=\frac{2^{n}(-i)^{n^{2}}}{n!}\omega^{n}.$$
This last one is given by the solution of Calabi conjecture (given by Aubin
in \cite{A2} and Yau in \cite{Y}).
\item
Finally, $L$ is a special Lagrangian manifold if we have
$i^{*}\Omega=dV_L$ (in the general case we have
$i^{*}\Omega=\lambda dV_L$ where $\lambda \in S^{1}$).
\end{itemize}
For a good initiation to the problem of Calabi conjecture and its
solution, one can refer to the book \cite{A1}. In our case and in
order to give examples of special Lagrangian manifolds, we use, as
it was done by R.L. Bryant \cite{B} and D. Joyce \cite{J}, a real
structure $c$ on $M$, i.e. an anti-holomorphic involution
$c:M\rightarrow M$ such that $c^{*}\omega=- \omega$ and
$c^{*}\Omega= \overline{\Omega }$. The set $L$ of fixed points of
this involution is a special Lagrangian manifold. $L$ can be viewed
as the real locus of $M$.

This paper gives an example of a $3$ dimensional complex
hypersurface of the Grassmannian $G_{2,4}\complex$ with vanishing
$C_1$. It is well known that its real locus is a special Lagrangian
manifold. Actually, it is a complete intersection of a quadric and a
quartic in $\proj_{5}\complex$. The first interpretation
(hypersurface of $G_{2,4}\complex$) does not require particular
skills in Algebraic Geometry. Therefore, both approaches will be
developed in this article. We also give an interpretation of the
obtained $3$-dimensional real locus as a 1-torus bundle over
$\proj_2\reel$. Using the second point of view (complete
intersection in $\proj_{5}\complex)$, we show that this is nothing
but $S^{3}/\mathbb{Z}_4$. The first section of this paper is devoted
to the description of the manifold $X$ (definition and calculus of
its first Chern class). In the second part, we stand and prove the
main result (theorem \ref{th}), giving both interpretations of the
real locus $L$ of $X$. Finally, we choose to describe two other
examples among a large class of manifolds which are similar to $X$.

\section{Definition and description of the manifold $X$}

Let us consider the following equation:
\begin{eqnarray}\label{eq1}
F(u,v)&=&(z_{0}z'_{1}-z_{1}z'_{0})^4+(z_{0}z'_{2}-z_{2}z'_{0})^4+(z_{0}z'_{3}-z_{3}z'_{0})^4\nonumber\\
&-&(z_{1}z'_{2}-z_{2}z'_{1})^4-(z_{1}z'_{3}-z_{3}z'_{1})^4-(z_{2}z'_{3}-z_{3}z'_{2})^4=0,
\end{eqnarray}
where $u=(z_0,z_1,z_2,z_3)$ and $u'=(z'_0,z'_1,z'_2,z'_3)$ are two
independent vectors of $\complex^4$. It is easy to see that this
equation depends only on the $2$-plane of $\complex^4$ given by $u$
and $u'$. Consequently, (\ref{eq1}) defines a subset of the
grassmannian $G_{2,4}\complex$, set of complex two dimensional
linear spaces of $\complex^{4}$.

\begin{lemma}\label{lem1}
The equation  (\ref{eq1}) defines a holomorphic hypersurface of
$G_{2,4}\complex$ which may be identified with a complete
intersection of a quadric and a quartic in $\proj_{5}\complex$.
\end{lemma}
\begin{proof}
In order to find the rank of the linear tangent map of $F$ in a
$G_{2,4}\complex$ classical coordinates system
$(\zeta_{i})_{i\in\{1,..,4\}}$, we consider~:

\begin{eqnarray}\label{eq2}
\left\{
\begin{array}{rr}
\zeta _{1}^{3}-\zeta _{4}\left( \zeta _{1}\zeta _{4}-\zeta _{2}\zeta
_{3}\right) ^{3}=0 &  \\
\zeta _{2}^{3}+\zeta _{3}\left( \zeta _{1}\zeta _{4}-\zeta _{2}\zeta
_{3}\right) ^{3}=0 &   \\
-\zeta _{3}^{3}+\zeta _{2}\left( \zeta _{1}\zeta _{4}-\zeta
_{2}\zeta
_{3}\right) ^{3}=0 &  \\
-\zeta _{4}^{3}-\zeta _{1}\left( \zeta _{1}\zeta _{4}-\zeta
_{2}\zeta _{3}\right) ^{3}=0 &
\end{array}
\right.
\end{eqnarray}
Using equation (\ref{eq1}), we firstly prove that the solutions of
this system are necessarily of norm one. Taking into account all
different cases, explicit computations lead to a contradiction
with the system (\ref{eq2}).\\
Let us give another description of $X$. To this end, we use the
classical identification of $G_{2,4}\complex$ with the quadric in
$\proj_{5}\complex$ given by~:
$$\eta_0\eta_5-\eta_1\eta_4+\eta_2\eta_3 =0,$$ via the map which
associates to every point
\begin{eqnarray}
\left(
\begin{array}{cc}
z_0 & z'_0\\
z_1 & z'_1\\
z_2 & z'_2\\
z_3 & z'_3
\end{array}
\right) \in G_{2,4}\complex \nonumber
\end{eqnarray}
described in the above coordinates, the point of
$\proj_{5}\complex$~:
\begin{eqnarray}
[\eta_{0}=z_{0}z'_{1}-z_{1}z'_{0},\eta_{1}=z_{0}z'_{2}-z_{2}z'_{0},\eta_{2}=z_{0}z'_{3}-z_{3}z'_{0},
\nonumber\\
\eta_{3}=z_{1}z'_{2}-z_{2}z'_{1},
\eta_{4}=z_{1}z'_{3}-z_{3}z'_{1},\eta_{5}=z_{2}z'_{3}-z_{3}z'_{2}].\nonumber
\end{eqnarray}
So, $X$ appears as a complete intersection of a quadric and a
quartic in $\proj_{5}\complex$, given by~:
\begin{eqnarray}\label{eq3}
\left\{
\begin{array}{rr}
\eta _{0}\eta _{5}-\eta _{1}\eta _{4}+\eta _{2}\eta _{3}=0  \\
\eta _{0}^{4}+\eta _{1}^{4}+\eta _{2}^{4}-\eta _{3}^{4}-\eta
_{4}^{4}-\eta _{5}^{4}=0
\end{array}
\right .
\end{eqnarray}
\end{proof}

\begin{lemma}\label{lem2}
$X$ is a manifold with vanishing first Chern class.
\end{lemma}

\begin{proof}

In a first time, we shall give a "self-contained" proof, using the
description  of $X$, given by equation (\ref{eq1}). Then, in a
second time, taking into account equation (\ref{eq3}), we give a
shorter proof which uses some concepts of Algebraic Geometry.

1) Let us prove that the determinant bundle $\Lambda^3T^*X$ is
trivializable by giving a $(3,0)$-holomorphic volume form $\Omega$
from a $(4,0)$-meromorphic form on $G_{2,4}\complex$, using
Poincar\é residue. It corresponds to the generalization of the
classical Cauchy residue at a point of a domain of $\complex$ to the
concept of residue in a hypersurface of a $n$ dimensional complex
manifold. If $\eta$ is a $n$-meromorphic form of such a manifold
which has first order poles on the hypersurface  $X$ locally defined
by the equation $g=0$, then $\eta$ may be locally written as:
$$\eta=\frac{\gamma \wedge dg}{g}+\delta ,$$
where $\gamma$ and $\delta$ are respectively $(n-1)$ and $n$ holomorphic
forms. Then, the restriction of $\gamma$ to $X$ is well defined as a
$(n-1)$ holomorphic form on $X$. We say that $\gamma$ is the
Poincar\'e residue of $\eta$.

In our case, let us consider the open chart set $U_{01}$ of
$G_{2,4}\complex$ where the $\{U_{ij} \}_{0\leq i<j\leq 3}$ are
given by~:
$$U_{ij}=\{(u,v)\in \complex^4 \times \complex^4 : u=(z_0,z_1,z_2,z_3)
,u'=(z'_0,z'_1,z'_2,z'_3),z_{i}z'_{j}-z_{j}z'_{i}\neq 0 \}.$$

The chart maps $\varphi_{ij}:U_{ij} \rightarrow \complex^4 \sim
M_2(\complex)$, are defined, similarly to $\varphi_{01}$ in the
following manner~:

\begin{eqnarray}
\varphi_{01}(u,v)= \left(
\begin{array}{cc}
z_2 & z'_2\\
z_3 & z'_3
\end{array}
\right) \times \left(
\begin{array}{cc}
z_0 & z'_0\\
z_1 & z'_1
\end{array}
\right)^{-1}= \left(
\begin{array}{cc}
\zeta _3 & \zeta _1\\
\zeta _4 & \zeta _2
\end{array}
\right) \nonumber
\end{eqnarray}

The expression of $F$ in $U_{01}$ is given by
$$f_{01}(\zeta_{1},\zeta_{2},\zeta_{3},\zeta_{4})=
1+\zeta_{1}^4+\zeta_{2}^4-\zeta_{3}^4-\zeta_{4}^4-(\zeta_{1}\zeta_{4}-\zeta_{2}\zeta_{3})^4
,$$

In $U_{01}$, let us consider~:
$$\eta = \frac{d\zeta_{1}\wedge d\zeta_{2}\wedge d\zeta_{3}\wedge
d\zeta_{4}}{1+\zeta_{1}^{4}+\zeta_{2}^{4}-\zeta_{3}^{4}
-\zeta_{4}^{4}-(\zeta_{1}\zeta_{4}-\zeta_{2}\zeta_{3})^{4}}.$$
$\eta$ is the local expression in the open chart set $U_{01}$ of a
$n$ meromorphic form globally defined on $G_{2,4}\complex$ whose
poles are along $X$. Indeed, the power $4$ is the correct one. This
can easily be seen in proceeding to the change of charts. There are
two types of change of charts in $G_{2,4}$. The first ones (or the
easiest) are of the form:
$$(\zeta_1,\zeta_2,\zeta_3,\zeta_4)\longrightarrow (-\frac{\zeta_2}{\zeta_{1}},\frac{1}{\zeta_{1}},-\frac{D}{\zeta_{1}},
\frac{\zeta_3}{\zeta_{1}}),$$ where
$D=(\zeta_{1}\zeta_{4}-\zeta_{2}\zeta_{3})$. The determinant of such
a change of charts is equal to $1/\zeta_{1}^{4}$, and so $\eta$ may
be extended to the new chart. The second change of charts is:
$$(\zeta_1,\zeta_2,\zeta_3,\zeta_4)\longrightarrow (-\frac{\zeta_4}{D},\frac{\zeta_2}{D},\frac{\zeta_3}{D},
-\frac{\zeta_1}{D}).$$ Its determinant (much more difficult to
compute) is equal to $1/D^4$, so, we may extend the expression to
the new chart. This result has been established in the general case
of $G_{p,p+q}\complex$ par J. Grivaux in \cite{G}. Using a skill
calculation, he has found the Einstein constant used in the K\"ahler
potential of Grassmannians. The Poincar\'e residue of
$$\eta = \frac{d\zeta_{1}\wedge d\zeta_{2}\wedge d\zeta_{3}\wedge d\zeta_{4}}{1+\zeta_{1}^{4}+\zeta_{2}^{4}-\zeta_{3}^{4}
-\zeta_{4}^{4}-(\zeta_{1}\zeta_{4}-\zeta_{2}\zeta_{3})^{4}}$$ is
given, in a submersion chart $X$ obtained as a sub-chart of
$(U_{01},\varphi_{01})$ considering the condition $\partial
f_{01}/\partial\zeta_1\neq 0$, by the expression:
$$\gamma = \frac{-d\zeta_{2}\wedge d\zeta_{3}\wedge d\zeta_{4}}{\partial f_{01}/\partial\zeta_1},$$
which defines a holomorphic $3$ volume form on $X$. This proves that $C_1(X)=0$. \\
2) Using the description (\ref{eq3}) of $X$, we can directly
establish this last result. Actually, a similar proof to the
preceding shows that a hypersurface of degree $m+1$ of
$\proj_{m}\complex$ has necessarily a vanishing $C_1$. Taking into
account that the sum of the degrees of a quadric and a quartic is
equal to 6, and because our intersection is complete in
$\proj_5\complex$, we obtain the result thanks to the adjunction
formula.
\end{proof}

\section{Lagrangian submanifold}

\begin{theorem}\label{th}
The real locus $L$ of $X$ is a special Lagrangian submanifold
endowed with a circle bundle structure over $\proj_2\reel$ which may
be identified with $S^{3}/\mathbb{Z}_{4}$.
\end{theorem}

\begin{proof}
In order to show that $L$ is a  special Lagrangian submanifold, we
use the result given by R.L. Bryant \cite{B}: the real locus of a
trivial first Chern class manifold (that is to say the fixed
points of an anti-holomorphic involution), when it is non-empty,
is a special Lagrangian manifold which is, according to the lemmas
\ref{lem1} and \ref{lem2}, the case of the manifold $X$ (the
anti-holomorphic involution we use here is the classical
conjugation). $L=L'/\sim$, where $L'$ is the set of the $2$-planes
of $\reel^{4}$ seen as pairs of independent vectors
$({u},{u}')\in\reel^4\times\reel^4$,
${u}=(x_0,x_1,x_2,x_3),
{u}'=(x'_0,x'_1,x'_2,x'_3)$ (defining a
point of $G_{2,4}\reel$) such that:
\begin{eqnarray*}
\left\{
\begin{array}{rr}
(x_{0}x'_{1}-x_{1}x'_{0})^4+(x_{0}x'_{2}-x_{2}x'_{0})^4+(x_{0}x'_{3}-x_{3}x'_{0})^4=1 \,\,\,\,(E_1)\\
(x_{1}x'_{2}-x_{2}x'_{1})^4+(x_{1}x'_{3}-x_{3}x'_{1})^4+(x_{2}x'_{3}-x_{3}x'_{2})^4=1
\,\,\,\,(E_2)
\end{array}
\right.
\end{eqnarray*}
and $\sim$ is the equivalence relation defined by:
$({u},{u}') \sim
({v},{v}')\in\reel^4\times\reel^4,$
if and only if
${v}=a{u}+c{u}'$ and
${v}'=b{u}+d{u}'$
with $ad-bc=\pm 1$.

$(E_1)$ and $(E_2)$ come from the equation (\ref{eq1}) and a
given normalization. So $L'$ is a submanifold of $\reel^8$,
diffeomorphic to the pull-back of $(1,1)\in\reel^2$ by the
submersion
$$\psi : (\reel\times\reel^3)\times(\reel\times\reel^3)\rightarrow\reel^2$$
defined by
$$\psi ((\alpha,{u_0}),(\alpha',{u_0}'))=(\|{u_0}\wedge
{u_0}'\|^{2},\|\alpha'{u_0}-\alpha{u_0}'\|^{2}).$$
In the former description
${u}=(\alpha,{u_0})\in\reel^4$ where
$\alpha = x_0$ and ${u_0}=(x_1,x_2,x_3)\in \reel^3$
is the projection of ${u}=(x_0,x_1,x_2,x_3)$ on
$\reel^3 = \{(0,x,y,z)\in \reel^{4}\}$ (it is the same for the prime
items). So the expression ${u_0}\wedge
{u_0}'$ is intrinsic. The result is more
difficult to obtain for
$\|\alpha'{u_0}-\alpha{u_0}'\|^{2}$,
after quotienting by $\sim$. Indeed, if
$${v}=a{u}+c{u}'=(\beta,{v_0})
\mbox{ and }
{v}'=b{u}+d{u}'=(\beta',{v_0}')
\mbox{ with } ad-bc=\pm 1,$$ we have
$$\beta=a\alpha+c\alpha',\,\,\,
\,\,\,\beta'=b\alpha+d\alpha',\,\,\,
{v_0}=a{u_0}+c{u_0}'\mbox{
and }
{v_0}'=b{u_0}+d{u_0}'.$$
So
\begin{eqnarray*}\beta'{v_0}-\beta{v_0}'&=&(b\alpha+d\alpha')
(a{u_0}+c{u_0}')-(a\alpha+c\alpha')
(b{u_0}+d{u_0}')\\
&=&(ad-bc)(\alpha'{u_0}-\alpha{u_0}')\\
&=&\pm (\alpha'{u_0}-\alpha{u_0}')
\end{eqnarray*}

The interpretation of $L$ as a circle bundle over $\proj_2\reel$
may be realized as follows:
\begin{itemize}
\item $(E_2)$ determines the basis of this bundle and corresponds to the fact that
the vectorial product of the vectors projections
${u}=(x_0,x_1,x_2,x_3)$ and
${u}'=(x'_0,x'_1,x'_2,x'_3)$ on $\reel^3 =
\{(0,x,y,z)\in \reel^{4}\}$ is unitary for a certain norm.
According to $(E_2)$ the vectors
${u_0}=(0,x_1,x_2,x_3)$ and
${u_0}'=(0,x'_1,x'_2,x'_3)$ are linearly
independent. Then they define a point of the Grassmannian
$G_{2,3}\reel$, that is to say a point of $\proj_2\reel$ (given by
their vectorial product).
\item $(E_1)$ describes the fibres above $\proj_2\reel$. To any
pair $({w},{w}')$ of independent
vectors of $\reel^3$, corresponds the circle in the basis
$\{{w},{w}'\}$ given by the equation
 $\|\alpha'{w}-\alpha{w}'\|^{2}=1$.
\end{itemize}

Using the second interpretation of $X$ (intersection of a quadric
and a quartic), and the notations used in the equations
(\ref{eq3}), $L$ may be identified with $S^{3}/\mathbb{Z}_{4}$ in
the following way~: It is the set of the vectors pairs of
$\reel^{3}$ $({u},{v})$ where ${u}=(\eta _{0},\eta _{1},\eta
_{2})$ and ${v}=(\eta _{5},-\eta _{4},\eta _{3})$ with vanishing
scalar product fulfilling the normalization condition
$$\eta _{0}^{4}+\eta _{1}^{4}+\eta _{2}^{4}=\eta _{5}^{4}+(-\eta
_{4})^{4}+\eta _{3}^{4}=1,$$ which is topologically equivalent to
$\|{u}\|=\|{v}\|$. If one considers the unit vectors
${u}'={u}/\|{u}\|$ and ${v}'={v}/\|{v}\|$, the triple
$({u}',{v}',{u}'\wedge{v}')$ is a direct orthonormal basis of
$\reel^{3}$ which can be viewed as a  matrix in $SO(3)$. The one
to one correspondence
$$[\eta_{0},..,\eta_{5}]\longrightarrow \{({u}',{v}'),
(-{u}',-{v}')\},$$ defines a map from $L$ in
$SO(3)/\mathbb{Z}_{2}$, where $\mathbb{Z}_{2}$ is the subgroup
$\{Id,\sigma\}$ of $SO(3)$ generated by the matrix $\sigma =\left(
          \begin{array}{ccc}
            -1 & 0 & 0 \\
           0 & -1 & 0\\
            0 & 0 & 1 \\
          \end{array}
        \right)$,
which corresponds to the map
$$({u}',{v}',{u}'\wedge{v}')
\longrightarrow
(-{u}',-{v}',{u}'\wedge{v}').$$
This map being a diffeomorphism, $L$ is diffeormorphic to
$SO(3)/\mathbb{Z}_{2}$ or equivalently $SU(2)/\mathbb{Z}_{4}$ that
is $S^{3}/\mathbb{Z}_{4}$, using the double cover
$SU(2)\longrightarrow SO(3)$. We can recover the description of $L$
(given in 1)) as a circle bundle over $\proj_{2}\reel$ by projecting
$$\begin{array}{ccc}
      SO(3)/\mathbb{Z}_{2} & \longrightarrow  & \proj_{2}\reel \\
      \{({u}',{v}',{u}'\wedge{v}'),
      (-{u}',-{v}',{u}'\wedge{v}')\}
      & \longrightarrow & \{{u}',-{u}'\}. \\
    \end{array}$$

\end{proof}

\section{Other examples}
To get other natural examples of special Lagrangian submanifolds
similar to $X$ it suffices to multiply, for example,
some factors by constants in equation (\ref{eq1}), keeping at least
one minus sign opposite to the others (in order to get a non empty
real locus). The nature of the obtained special Lagrangian
submanifolds may be quite different of the submanifold $L$.\\
1) Let us consider the following example whose
interpretation is of similar interest to the above one:
\begin{eqnarray}\label{eq7}
(z_{0}z'_{1}-z_{1}z'_{0})^4-(z_{0}z'_{2}-z_{2}z'_{0})^4-(z_{0}z'_{3}-z_{3}z'_{0})^4
-(z_{1}z'_{2}-z_{2}z'_{1})^4\nonumber\\
-(z_{1}z'_{3}-z_{3}z'_{1})^4-2(z_{2}z'_{3}-z_{3}z'_{2})^4=0.
\end{eqnarray}
where $u=(z_0,z_1,z_2,z_3)$ and $u'=(z'_0,z'_1,z'_2,z'_3)$ are two
independent vectors of $\complex^4$, defining a $2$-plane of
$\complex^{4}$, that is to say a point of $G_{2,4}\complex$. As it
was realized above for $X$, we prove that this hypersurface is a
trivial first Chern class submanifold. Its real locus which is a
special Lagrangian submanifold is given (using a good
normalization) by the equations
$$\hat{E}_1:\,\,\,(x_{0}x'_{1}-x_{1}x'_{0})^4=1$$
and
\begin{eqnarray*}
\hat{E}_2:\,\,\,\,(x_{1}x'_{2}-x_{2}x'_{1})^4+(x_{1}x'_{3}-x_{3}x'_{1})^4+2(x_{2}x'_{3}-x_{3}x'_{2})^4\\
+(x_{0}x'_{2}-x_{2}x'_{0})^4+(x_{0}x'_{3}-x_{3}x'_{0})^4=1,
\end{eqnarray*}
quotiented by the action $\sim$ described above for $X$.
\begin{itemize}
\item $\hat{E}_1$ indicates that $L'$ is a subset of the open
chart set $(x_{0}x'_{1}-x_{1}x'_{0})\neq 0$ of $G_{2,4}\reel$.
\item$\hat{E}_2$ may be written as:
\begin{eqnarray*}
\hat{E}_2&:&(x_{1}x'_{2}-x_{2}x'_{1})^4+(x_{1}x'_{3}-x_{3}x'_{1})^4+(x_{2}x'_{3}-x_{3}x'_{2})^4\\
&&(x_{2}x'_{3}-x_{3}x'_{2})^4+(x_{0}x'_{2}-x_{2}x'_{0})^4+(x_{0}x'_{3}-x_{3}x'_{0})^4=1,
\end{eqnarray*}
\end{itemize}
which is a set diffeomorphic to
$$\|{u}_0\wedge{u}'_0\|^{2}+\|{u}_1\wedge
{u}'_1\|^{2}=1,$$ where ${u}_0$ and
${u}_1$ are the projections of
${u}=(x_0,x_1,x_2,x_3)\in \reel^4$ respectively on
 $\{x_0=0\}$ and $\{x_1=0\}$ seen as
$\reel^3$ (we have the same situation for the prime items).

2) A second example is given for any $u=(z_0,z_1,z_2,z_3)$ and
$u'=(z'_0,z'_1,z'_2,z'_3)$, independent vectors of $\complex^4$
(seen as a point of $G_{2,4}\complex$), by:
\begin{eqnarray}\label{eq8}
(z_{0}z'_{1}-z_{1}z'_{0})^4+(z_{0}z'_{2}-z_{2}z'_{0})^4-(z_{0}z'_{3}-z_{3}z'_{0})^4
-2(z_{1}z'_{2}-z_{2}z'_{1})^4\nonumber\\
-2(z_{1}z'_{3}-z_{3}z'_{1})^4-2(z_{2}z'_{3}-z_{3}z'_{2})^4=0.
\end{eqnarray}
and defining a hypersurface of $G_{2,4}\complex$ with trivial
first Chern class. Its real locus $R$, which is a special
Lagrangian manifold is given, after the normalization used in both
former examples, by the equations:
\begin{eqnarray*}
\left\{
\begin{array}{ll}
(x_{0}x'_{1}-x_{1}x'_{0})^4+(x_{0}x'_{2}-x_{2}x'_{0})^4=1 \,\,\,\,(\tilde{E}_1)\\
(x_{0}x'_{3}-x_{3}x'_{0})^4+2(x_{1}x'_{2}-x_{2}x'_{1})^4
+2(x_{1}x'_{3}-x_{3}x'_{1})^4+2(x_{2}x'_{3}-x_{3}x'_{2})^4=1
\,\,\,\,(\tilde{E}_2)
\end{array}
\right.
\end{eqnarray*}
again quotiented by $\sim$ (as described above). Using a
normalization, this system is equivalent to the following one:
\begin{eqnarray*}
\left\{
\begin{array}{ll}
(x_{0}x'_{1}-x_{1}x'_{0})^4+(x_{0}x'_{2}-x_{2}x'_{0})^4+(x_{0}x'_{3}-x_{3}x'_{0})^4=1 \,\,\,\,(\tilde{E'}_1)\\
2(x_{0}x'_{3}-x_{3}x'_{0})^4+2(x_{1}x'_{2}-x_{2}x'_{1})^4
+2(x_{1}x'_{3}-x_{3}x'_{1})^4+2(x_{2}x'_{3}-x_{3}x'_{2})^4=1
\,\,\,\,(\tilde{E'}_2)
\end{array}
\right.
\end{eqnarray*}
As it was noticed in the case of $X$, $(\tilde{E'}_1)$ describes
the circle
$\|x_{0}{u}'_0-x'_{0}{u}_0\|^{2}=1$.
Recall that ${u}_0$ is the projection of
${u}=(x_0,x_1,x_2,x_3)\in \reel^4$ on $\{x_0=0\}$
(idem for the prime items). The interpretation of $(\tilde{E'}_2)$
using classical geometric objects seems to be a little more
difficult. However on can affirm that the set described by
$(\tilde{E'}_2)$ contains the ball
$\|{u}_0\wedge{u}'_0\|^{2}<1/2.$


\end{document}